# Construction of sphere maps with given degrees and a new proof of Morse index formula

Xiao-Song Yang

Department of Mathematics,
Huazhong University of Science and Technology.
Wuhan, 430074, China

***Abstract***
This note presents a procedure of constructing a higher dimensional sphere map from a lower dimensional one and gives an explicit formula for smooth sphere map with a given degree. As an application a new proof of a generalized Poincare-Hopf theorem called Morse index formula is also presented.

## 1 Introduction

The degree of (continuous) smooth map $f: S^n \to S^n$ is a well established topic appeared in many text books and a vast literature. Nonetheless, how to construct specific maps from $S^n$ to $S^n$ with given dimension and topological degree number seems seldom appear in the literature, although it is interesting and meaningful. In this note we present a procedure of constructing a higher dimensional sphere map from a lower dimensional one and give an explicit formula for sphere map with given degree number. Thus for each homotopic class of continuous map $f: S^n \to S^n$ we can have an explicit representative map. As an application, we give a new proof of the generalized Poincare-Hopf theorem on indices of continuous vector fields, which was obtained by Morse [2], and rediscovered by Pugh [4] and Gottlleb [1]. For convenience, this generalized theorem is called Morse index formula following Grottleb [1].

## 2 A fundamental lemma

In this section we prove a lemma that is essential to our construction of higher dimensional sphere map by lower dimensional sphere map and important for the new proof the Morse index formula.

Let $V(x)$ be a continuous vector field defined on $U \subset R^{n+1}$, where $U$ is an open neighborhood of $x=0$. It is assumed that $x=0$ is the only zero of $V(x)$ on $U$. As usual, one can define the index $ind(V,0)$ of $V(x)$ at $x=0$ to be the degree of the map $f_V : S_r^n \to S^n$:

$$f_V(x) = \frac{V(x)}{\|V(x)\|}, \quad x \in S_r^n$$

where $S_r^n \subset U$ is a sphere centered at $x=0$ with radius $r$.

Let $R_0^n = \{x = (x_1,...,x_n,x_{n+1}) \in R^{n+1} : x_{n+1} = 0\}$ and suppose that the set $U^n = U \bigcap R_0^n$ is invariant for $V(x)$, i.e., $V(x)$ is tangent to $U^n$ at every point $x \in U^n$, which is equivalent to $V_{n+1}(x) = 0$ for every $x \in U^n$ when $V(x)$ is of the form

$$V(x) = (V_1(x),...,V_n(x),V_{n+1}(x))$$

It is apparent that $\bar{V} = V|U^n$, the restriction of $V(x)$ to $U^n$, is a tangent vector field on $U^n$, which also induces a lower dimensional sphere map

$$\bar{f}_{\bar{V}} : S_r^{n-1} \to S^{n-1} = S^n \bigcap R_0^n$$

$$\bar{f}_{\bar{V}}(x) = \frac{\bar{V}(x)}{\|\bar{V}(x)\|}, \quad x \in S_r^{n-1} = S_r^n \bigcap U^n .$$

The index $ind(\bar{V},0)$ of $\bar{V}(x)$ is defined to be the degree of the above map.

Now we have the following observations.

**Lemma 2.1** Let $V(x)$ be a continuous vector field defined on $U \subset R^{n+1}$. Suppose that $V(x)$ is tangent to $U^n$, Then

(i) $ind(V,0) = -ind(\bar{V},0)$, if $x_{n+1}V_{n+1}(x) < 0$ for $x_{n+1} \neq 0$

(ii) $ind(V,0) = ind(\bar{V},0)$, if $x_{n+1}V_{n+1}(x) > 0$ for $x_{n+1} \neq 0$

***Proof*** First we consider the case (i). Let

$$W(x) = (V_1(x),...,V_n(x),-x_{n+1}), \quad x \in U$$

Without loss of generality we can assume that the first $n$ components $V_1(x),...,V_n(x)$ are smooth (or one can smooth them), so that $\bar{f}_{\bar{V}}$ is smooth.

First we show that

$$ind(V,0) = ind(W,0).$$

To see this, construct a homotopy

$$H(x,t): S^n_{r^n} \times [0,1] \to S^n$$

$$H(x,t) = \frac{(1-t)V + tW}{\|(1-t)V + tW\|}$$

From the identity

$$(1-t)V + tW = (V_1(x),\ldots,V_n(x), -(1-t)x_{n+1} + V_{n+1}(x))$$

it is easy to see that this homotopy makes sense as long as $r$ is sufficiently small.

To see what the index $ind(W,0)$ is, let us consider the degree of $f_W : S^n_r \to S^n$.

$$f_W(x) = \frac{W(x)}{\|W(x)\|}, \quad x \in S^n_r$$

Since $U^n$ is also invariant for $W(x)$, one has

$$\bar{f}_{\bar{W}} : S^{n-1}_r \to S^{n-1} = S^n \cap R^n_0$$

$$\bar{f}_{\bar{W}}(x) = \frac{\bar{W}(x)}{\|\bar{W}(x)\|}, \quad x \in S^{n-1}_r = S^n_r \cap U^n.$$

where $\bar{W}(x)$ is the restriction of $W(x)$ to $U^n$, it is just $\bar{V} = V|U^n$. Therefore

$$\bar{f}_{\bar{W}} = \bar{f}_{\bar{V}}.$$

Let $p \in S^{n-1} = S^n \cap R^n_0$ be a regular value of $\bar{f}_{\bar{V}}$, then by definition

$$\deg(\bar{f}_{\bar{V}}, p) = \sum_{y \in \bar{f}_{\bar{V}}^{-1}(p)} \operatorname{sgn} D\bar{f}_{\bar{V}}(y)$$

where sgn $D\bar{f}_{\bar{V}}(y)$ = sign of $\det(D\phi_y \circ \bar{f}_{\bar{V}} \circ \varphi^{-1}_{p\,p}(\varphi_p(y)))$, and $(\phi_y, U_y)$ is a chart around $x \in S^{n-1}_r$ while $(\varphi_p, U_p)$ is a chart around $p \in S^{n-1}$.

For each $y \in \bar{f}_{\bar{V}}^{-1}(p)$, let

$$\bar{y} = (y,0) \in S^n_r \ \text{ and } \ \bar{p}(p,0) \in S^n$$

Then it is easy to extend the chart $(\phi_y, U_y)$ to a chart $(\phi_{\bar{y}}, U_y \times (-1,1))$ and $(\varphi_p, U_p)$ to $(\varphi_{\bar p}, U_p \times (-1,1))$ with $(\phi_{\bar y}|U_y \times \{0\} = \phi_y$ and $(\varphi_{\bar p}|U_p \times \{0\} = \varphi_p$.

Now

$$Df_W(\bar y) = D\phi_{\bar y} \circ f_W \circ \varphi_{\bar p}^{-1}(\varphi_{\bar p}(\bar y))$$

$$= \begin{bmatrix} D\phi_y \circ \bar f_{\bar V} \circ \varphi^{-1}_{p\,p}(\varphi_p(y)) & 0 \\ 0 & -1 \end{bmatrix}$$

Keeping in mind that $W(x)$ is tangent to $U^n$ we have

$$f_W^{-1}(\bar p) = \bar f_{\bar V}^{-1}(p).$$

Therefore

$$\deg(f_W, \bar p) = \sum_{\bar y \in f_W^{-1}(\bar p)} \operatorname{sgn} Df_W(\bar y)$$

$$= -\sum_{y \in \bar f_{\bar V}^{-1}(p)} \operatorname{sgn} D\bar f_{\bar V}(y) = -\deg(\bar f_{\bar V}, p)$$

Consequently $ind(V,0) = ind(W,0) = -ind(\bar V,0)$, thus we prove the first case.

For case (ii), we can in the same manner define a new vector field

$$W(x) = (V_1(x), \ldots, V_n(x), x_{n+1}), \quad x \in U$$

to show the identity

$$\deg(f_W, \bar p) = \sum_{\bar y \in f_W^{-1}(\bar p)} \operatorname{sgn} Df_W(\bar y)$$

$$= \sum_{y \in \bar f_{\bar V}^{-1}(p)} \operatorname{sgn} D\bar f_{\bar V}(y) = \deg(\bar f_{\bar V}, p) \blacksquare$$

By virtue of **Lemma 2.1** it is easy to propose a procedure of constructing an explicit smooth map from $n-$dimensional sphere to itself with a given topological degree, as shown in the next section..

## 3 Constructing $m-$degree sphere maps

The purpose of this section is to present an explicit smooth map from $S^n$ to $S^n$ with given topological degree. First we consider the construction of $m-$degree smooth map from $S^1$ to $S^1$.

For an integer $m \geq 0$, let

$$P_m(x_1, x_2) = \operatorname{Re}(x_1 + ix_2)^m, \qquad Q_m(x_1, x_2) = \operatorname{Im}(x_1 + ix_2)^m$$

Consider the planar vector field

$$V_m(x_1, x_2) = (P_m(x_1, x_2), Q_m(x_1, x_2)),$$

which has the origin as its only zero. With this vector field define the following map

$$\alpha_1^m : S^1 \to S^1$$

$$\alpha_1^m = \frac{V_m(x_1, x_2)}{\|V_m(x_1, x_2)\|}, \qquad (x_1, x_2) \in S^1 \subset R^2.$$

It is easy to prove the following fact

**Lemma 3.1** The map $\alpha_1^m : S^1 \to S^1$ is smooth and of $m-$degree.

Now consider the vector field $\overline{V}_m(x)$ defined on $R^{n+1}$:

$$\overline{V}_m(x) = (x_1, ..., x_{n-1}, P_m(x_n, x_{n+1}), Q_m(x_n, x_{n+1})), \quad x = (x_1, ..., x_{n-1}, x_n, x_{n+1}) \in R^{n+1}$$

and define the map $\alpha_n^m : S^n \to S^n$:

$$\alpha_n^m = \frac{\overline{V}_m(x)}{\|\overline{V}_m(x)\|}, \qquad x = (x_1, ..., x_{n-1}, x_n, x_{n+1}) \in S^n \subset R^{n+1}.$$

We have the following statement

**Proposition 3.1** Every continuous map $f : S^n \to S^n$ with degree $m \geq 0$ is homotoppic to

$$\alpha_n^m : S^n \to S^n.$$

Now consider the following vector field $\overline{V}_{-m}(x)$ defined on $R^{n+1}$

$$\overline{V}_{-m}(x) = (-x_1, ..., x_{n-1}, P_m(x_n, x_{n+1}), Q_m(x_n, x_{n+1})), \quad x = (x_1, ..., x_{n-1}, x_n, x_{n+1}) \in R^{n+1}$$

and define the map $\alpha_n^m : S^n \to S^n$:

$$\alpha_n^{-m} = \frac{\overline{V}_{-m}(x)}{\|\overline{V}_{-m}(x)\|}, \qquad x = (x_1, ..., x_{n-1}, x_n, x_{n+1}) \in S^n \subset R^{n+1}$$

We have the following statement

**Proposition 3.2** Every continuous map $f : S^n \to S^n$ with degree $-m \leq 0$ is homotoppic to

$$\alpha_n^{-m} : S^n \to S^n.$$

Keeping in mind the **Lemma 3.1** and the Hopf theorem.[2], the above two propositions can be easily proved by induction.

### 4 A new proof of the Morse index formula

Let $M$ be a compact smooth manifold with or without boundary $\partial M$, and let $V(x)$ be a continuous vector field on $M$ with only isolated zeros. If $\partial M$ is nonempty then $V(x)$ is assumed point outward at points of $\partial M$. Under these assumptions there holds the classical Poincare-Hopf index theorem, i.e., the index sum of the zeros of $V(x)$ equals to the Eular characteristic of the manifold $M$.

M. Morse [2] gave a more general index formula of a vector field $V(x)$ on $M$ with nonempty boundary $\partial M$ under a generic condition.

Suppose that $V(x)$ is a vector field on $M$ with zeros in $M - \partial M$, and all the zeros are isolated. Let $\partial_M$ be the component of $\partial M$ where $V(x)$ points inward, and $\partial_V(x)$ be the tangential component of $V(x)$ to $\partial M$ when restricted to $\partial_M$. The vector field $\partial_V(x)$ is assumed have only isolated zeros (if any) in $\partial_M$.

Denote by $Ind(V)$ the index sum of zeros of $V(x)$ on $M$ and $Ind(\partial_V)$ the index sum of zeros of $\partial_V(x)$ on $\partial_M$, then under the above assumptions, the so called Morse index formula can be stated as follows [1].

$$Ind(V) + Ind(\partial_V) = \chi(M)$$

In this section we present a new proof based on the lemma given in Section 2. This new proof seems more understandable and is more general in the sense that it does not needs the assumption that zeros are non-degenerate, which is essential to the arguments by Pugh [4].

***Proof*** of the Morse index formula

Since $V(x)$ has no zero on $\partial M$, there is collar neighborhood $\alpha \subset M$ of $\partial M$ in which $V(x)$ has no zero. Also there is a diffeomorphism

$$\phi : \alpha \to \partial M \times [0,1)$$

in view of the collar theorem.

Now $\phi_* V$ can be written as

$$\phi_* V = (V_t, V_v)$$

Where $V_t$ is tangent to $\partial M$ when restricted to $\partial M$, and $V_v(y) > 0$ or $V_v(y) < 0$ for

$$y = (\bar{y}, s) \in \partial M \times [0,1)$$

Consider a $C^{\infty}$ function

$$h : [0,1) \to [0,1]$$

satisfying

$$h(0) = 0,\ h(s) > 0,\ s \in (0, \tfrac{1}{2}],\ h(s) = 1,\ s \in [\tfrac{1}{2}, 1).$$

By means of this function one can construct a new vector field on $M$ as follows

$$\tilde{V}(x) = V(x), \quad \forall x \in M - U_c$$

$$\tilde{V}(x) = \phi_*^{-1}(V_t, h(s)V_v)(x) \ , \ x = \phi^{-1}(\bar{y}, s)$$

Clearly $\tilde{V}$ is smooth and is tangent to $\partial M$ at points in $\partial M$. In addition, $\tilde{V}$ has the same isolated zeros as $V$.

Now using the routing technique as did in the literature, we double the manifold $M$ and the corresponding vector field $\tilde{V}$, yielding

$$\overline{M} = M \bigcup M_c,\ \overline{V} = \tilde{V} \cup \tilde{V}_c$$

where $M_c$ is the copy of $M$ and $\tilde{V}_c$ is the copy of $\tilde{V}$ on $M_c$. In addition $\partial M = M \bigcap M_c$, and $\overline{V}$ is tangent to $\partial M$.

Now the collar embedding $\phi : U_c \to \partial M \times [0,1)$ yields

$$\bar{\phi} : \overline{U} \to \partial M \times (-1,1),\ \overline{U} = \alpha \cup \alpha_c,\ \alpha_c \text{ is the copy of } \alpha$$

Since by assumption $\overline{V}\big|\partial M$ has isolated zeros (if any) as a vector field on $\partial M$, these isolated zeros are also isolated zeros of $\overline{V}$ on $M$. Because the index of an isolated zero is a local property, it is convenient to consider $\overline{V}$ ( in fact $\bar{\phi}_*\overline{V}$) with only zero located at the $(0,0)$ of $U^{n-1} \times (-1,1)$, where $U^{n-1} \subset R^{n-1}$ is an open neighborhood of the origin.

Let $y = (\bar{y}, s) \in U^{n-1} \times (-1,1)$, then $\bar{\phi}_*\overline{V}$) can be written as

$$\bar{\phi}_*\overline{V}(y) = \begin{cases} \{V_t(\bar{y}, s), V_v(\bar{y}, s)) & s \geq 0 \\ \{V_t(\bar{y}, -s), -V_v(\bar{y}, -s)), & s \leq 0 \end{cases}$$

with $V_v(\bar{y},0)) = 0$.

By lemma 2.1, it is easy to see that

$$ind(\bar{\phi}_*\bar{V},0) = ind(\bar{\phi}_*\bar{V}\big|U^{n-1},0) \quad \text{if} \quad sV_v(\bar{y},s) > 0 \quad \text{for} \quad s \neq 0,$$

and

$$ind(\bar{\phi}_*\bar{V},0) = -ind(\bar{\phi}_*\bar{V}\big|U^{n-1},0) \quad \text{if} \quad sV_v(\bar{y},s) < 0 \quad \text{for} \quad s \neq 0$$

Therefore

$$ind(\bar{V},p) = ind(\bar{V}\big|\partial M,p) \quad \text{if} \quad p \in \partial_+ M$$

and

$$ind(\bar{V},p) = -ind(\bar{V}\big|\partial M,p) \quad \text{if} \quad p \in \partial_- M$$

where $\partial_+ M$ is the component of $\partial M$ where $V(x)$ points outward. It follows from the well known Poincare-Hopf theorem that

$$\begin{aligned} Ind\bar{V} &= \sum_{p\in M} ind(\tilde{V},p) + \sum_{p\in M_c} ind(\tilde{V}_c,p) + \sum_{p\in\partial M} ind(\bar{V},p) \\ &= \chi(\bar{M}) = 2\chi(M) - \chi(\partial M) \end{aligned}$$

Since

$$\sum_{p\in\partial M} ind(\bar{V},p) = -\sum_{p\in\partial_+ M} ind(\bar{V}\big|\partial M,p) + \sum_{p\in\partial_- M} ind(\bar{V}\big|\partial M,p)$$

and

$$ind(\tilde{V},p) = ind(V,p).$$

One has

$$\begin{aligned} Ind\bar{V} &= 2\sum_{p\in M} ind(V,p) + 2\sum_{p\in\partial_- M} ind(\bar{V}\big|\partial M,p) \\ &\quad -(\sum_{p\in\partial_- M} ind(\bar{V}\big|\partial M,p) + \sum_{p\in\partial_+ M} ind(\bar{V}\big|\partial M,p) \end{aligned}$$

Note that $\bar{V}\big|\partial_- M = \partial_- V$ and

$$\sum_{p\in\partial_- M} ind(\bar{V}\big|\partial M,p) + \sum_{p\in\partial_+ M} ind(\bar{V}\big|\partial M,p) = \sum_{p\in\partial M} ind(\bar{V}\big|\partial M,p) = \chi(\partial M)$$

Therefore

$$\begin{aligned} &2\sum_{p\in M} ind(V,p) + 2\sum_{p\in\partial_- M} ind(\bar{V}\big|\partial M,p) - \chi(\partial M) \\ &= Ind\bar{V} = 2\chi(M) - \chi(\partial M) \end{aligned}$$

It follows that

$$\sum_{p\in M} ind(V,p) + \sum_{p\in\partial_M} ind(\bar{V}\big|\partial M, p) = \chi(M)$$

That is

$$Ind(V) + Ind(\partial_V) = \chi(M)$$

This completes the proof of Morse index lemma.■

**Acknowledgements** This work is supported in part by National Natural Science Foundation of China （10972082）.